\documentclass[11pt,reqno]{amsart}
\usepackage[T1]{fontenc}
\usepackage{lmodern}
\usepackage{amsmath,amssymb,amsthm,mathtools,microtype}
\usepackage{cite,etoolbox}
\usepackage{xcolor}
\usepackage[colorlinks=true,linkcolor=blue!45!black,citecolor=blue!45!black,urlcolor=blue!45!black]{hyperref}
\hypersetup{pdftitle={A note on the cube inequality},pdfauthor={Yuchen Bi}}
\numberwithin{equation}{section}
\newtheorem{theorem}{Theorem}[section]
\newtheorem{proposition}[theorem]{Proposition}
\newtheorem{lemma}[theorem]{Lemma}

\theoremstyle{definition}
\newtheorem{definition}[theorem]{Definition}
\theoremstyle{remark}
\newtheorem{remark}[theorem]{Remark}
\title[A note on the cube inequality]{A note on the cube inequality}
\author{Yuchen Bi}
\address{Mathematisches Institut, Universit\"at Freiburg, Germany}
\email{yuchen.bi@math.uni-freiburg.de}
\date{}
\allowdisplaybreaks
\makeatletter
\patchcmd{\@maketitle}{\dimen@34\p@}{\dimen@20\p@}{}{\errmessage{Title spacing patch failed}}
\makeatother
\begin{document}
\begin{abstract}
In this note, we adapt ideas from No.~336 of OpenAI's October~6, 2026 mathematics release to prove the weighted cube inequality in all dimensions.
\end{abstract}

\maketitle

\section{Introduction}

The cube inequality was formulated by Gromov~\cite[Section~3.8]{Gromov}. Let \(I=[-1,1]\) and let \((X^n,\partial X,g)\) be a compact, smooth and oriented Riemannian manifold. Assume there is a continuous map
\[
F=(F_1,\ldots,F_n):(X,\partial X)\longrightarrow (I^n,\partial I^n)
\]
with nonzero degree. Write
\[
B_i^\pm=\{F_i=\pm 1\},\qquad d_i=\operatorname{dist}_g(B_i^-,B_i^+).
\]
Then
\begin{equation}\label{eq:classicalcube}
R_g\ge\sigma>0\Longrightarrow\sigma\le 4\pi^2\frac{n-1}{n}\sum_{i=1}^nd_i^{-2}.
\end{equation}
Gromov proved the inequality for \(n\le 8\) using \(\mu\)-bubbles~\cite[Sections~5.2 and~5.4]{Gromov}, and conjectured it for \(n\ge 9\). Wang--Xie--Yu~\cite{WXY} resolved the conjecture for spin manifolds using Callias-type operators. For foundational work on spinorial methods in scalar curvature, see Gromov--Lawson~\cite{GL,GL83} and Witten~\cite{Witten}; see also Cecchini--Zeidler~\cite{CZ} for related band width estimates.

The weighted analogue of \eqref{eq:classicalcube} involves the weighted scalar curvature
\begin{equation}\label{eq:weightedscalar}
\mathcal S_N(g,f)=R_g+2\Delta_gf-\left(1+\frac{1}{N-n}\right)|df|_g^2,
\end{equation}
where \(f\in C^\infty(X)\) and \(N<0\) or \(n<N\le\infty\), with coefficients at \(N=\infty\) defined by limits. It asserts that \(\mathcal S_N(g,f)\ge\sigma>0\) implies
\begin{equation}\label{eq:cube}
\sigma\le 4\pi^2\frac{N-1}{N}\sum_{i=1}^nd_i^{-2}.
\end{equation}
For \(f=0\), letting \(N\downarrow n\) recovers \eqref{eq:classicalcube}.

Jintian Zhu and the author~\cite[Theorem~B]{BZ} proved \eqref{eq:cube} under the condition
\[
\frac{1}{N}\le\begin{cases}
\dfrac{1}{n},&n\le13,\\[5pt]
\dfrac{51-2n}{23n+26},&n\ge14.
\end{cases}
\]
For \(f=0\), this gives the sharp inequality \eqref{eq:classicalcube} for \(n\le 13\) and the nonsharp coefficient \(25(n-1)/(23n+26)\) for \(n\ge 14\), which suffices to show that closed enlargeable manifolds admit no metric of positive scalar curvature. The method of proof is a recent conformal blow-up descendant of the Schoen--Yau minimal hypersurface approach~\cite{SY,SYPMT}, first developed for the positive mass theorem~\cite{HSY,BHHSZ,BW}.

This note grew out of studying OpenAI's \emph{Positive scalar curvature and uniform codimension-two width}~\cite{OAI}. We prove the weighted cube inequality for the remaining values of \(N\) using another descendant of the Schoen--Yau approach, based on blow-up of the prescribed mean curvature. The treatment of singularities in Sections~3 and~4 largely follows \cite[Appendices~C--E]{OAI}.

\begin{theorem}\label{thm:cube}
Let \(N<0\) or \(n<N\le\infty\). If \(f\in C^\infty(X)\) satisfies \(\mathcal S_N(g,f)\ge\sigma>0\), then \eqref{eq:cube} holds. The coefficient \((N-1)/N\) is optimal.
\end{theorem}

We prove Theorem~\ref{thm:cube} when
\[
\frac{N-1}{N}<\frac{(7-2)^2}{4}.
\]
Outside this range, the theorem follows from \cite[Theorem~B]{BZ} with \(\lambda=1/N\).

The proof proceeds by dimension descent. An \(m\)-dimensional descent datum consists of a weighted manifold \((X_m,g,f_m)\) (where \(g\) denotes \(g|_{X_m}\) by abuse of notation) and a map
\[
(\Phi_1,\ldots,\Phi_m):(X_m,\partial X_m)\longrightarrow([-\infty,\infty]^m,\partial[-\infty,\infty]^m)
\]
of nonzero degree, satisfying
\begin{equation}\label{eq:intro-Phi}
2|d\Phi_j|\le\frac{N}{N-1}\Phi_j^2+\sigma_j\quad\text{a.e. on }\{|\Phi_j|<\infty\},\qquad j=1,\ldots,m,
\end{equation}
and
\begin{equation}\label{eq:intro-curvature}
\mathcal S_N(g,f_m)\ge s+\sum_{j=1}^m\sigma_j\quad\text{outside }\bigcup_{j=1}^m(U_j^+\cup U_j^-).
\end{equation}
Here \(s>0\) is a constant, the \(\sigma_j\) are smooth positive functions, and the relatively open sets \(U_j^\pm\) cover \(\partial X_m\), with \(\Phi_j=\pm\infty\) near \(\overline{U_j^\pm}\).

Using \(\Phi_m\), we construct a \(\mu\)-bubble \(\Sigma\) (as made precise in Lemma~\ref{lem:minimizer}) avoiding \(U_m^+\cup U_m^-\). For \(2\le m\le7\), \(\Sigma\) is smooth and the stability inequality yields the descent by changing the weight. For \(m\ge8\), we use the Hardy estimate in Section~\ref{sec:hardy} to modify \(\Phi_j\) to \(\widehat\Phi_j\), \(j<m\), so that \(\widehat\Phi_j=\pm\infty\) on open sets \(\widehat U_j^\pm\supset U_j^\pm\) whose union contains the singular set of \(\Sigma\) in \(X_m\), while \(\widehat\Phi_j\) satisfies \eqref{eq:intro-Phi} with smooth functions \(\widehat\sigma_j\ge\sigma_j\) in place of \(\sigma_j\). The sum \(\sum_{j<m}(\widehat\sigma_j-\sigma_j)\) can be absorbed in the stability inequality, yielding a smooth weight \(f_{m-1}\) on the hypersurface \(X_{m-1}\) obtained by smoothing \(\Sigma\) within \(\bigcup_{j<m}(\widehat U_j^+\cup\widehat U_j^-)\). This gives an \((m-1)\)-dimensional descent datum with functions \(\widehat\Phi_j|_{X_{m-1}}\), \(j<m\), and an arbitrarily small loss in the curvature lower bound. The degree is unchanged, since \(\widehat\Phi_j=\Phi_j\) near \(\partial X_m\) for \(j<m\) and \(X_{m-1}\) separates \(U_m^+\) from \(U_m^-\).

If \eqref{eq:cube} failed, \((X,g,f)\) together with a map \((\Phi_1,\ldots,\Phi_n)\) obtained by modifying the functions
\[
-\frac{2\pi(N-1)}{Nd_i}\cot\left(\frac{\pi}{d_i}\min\{\operatorname{dist}_g(x,B_i^-),d_i\}\right),\qquad i=1,\ldots,n,
\]
to equal \(\pm\infty\) near \(B_i^\pm\) would form an \(n\)-dimensional descent datum with \(s>0\). Iterating the descent would then give a contradiction in dimension one.

\subsection*{AI disclosure}
This note arose from several discussions with \mbox{ChatGPT} (GPT-6 Astra) about \cite{OAI}. \mbox{ChatGPT} assisted in drafting and revising the manuscript and in finding references. The author takes full responsibility for this paper.

\section{Stability and change of weight}

We derive the stability inequality for \(\mu\)-bubbles.

Let \((M^m,g)\) be a Riemannian manifold with \(m\ge2\), and let \(f\in C^\infty(M)\). Fix \(N<0\) or \(m<N\le\infty\). Write
\[
d\mathfrak m=e^{-f}d\mu_g,\qquad \Delta_f=\Delta_g-\langle df,d \cdot \rangle_g.
\]

For \(\Phi\in C^\infty(M)\), let \(\Omega\) be a set of locally finite perimeter that minimizes
\begin{equation}\label{eq:functional}
\mathcal J(E;U)=\int_{\partial^*E\cap U}d\mathfrak m-\int_{E\cap U}\Phi d\mathfrak m
\end{equation}
among sets \(E\) of locally finite perimeter with \(E\triangle\Omega\Subset U\), for every open set \(U\Subset M\). Set \(\Sigma=\operatorname{spt}|D\chi_\Omega|\). On its regular part \(\Sigma_{\mathrm{reg}}\), let \(\nu\) be the outward unit normal, \(A\) the second fundamental form, and \(H=\operatorname{tr}A=\operatorname{div}_\Sigma\nu\). We also write \(d\mathfrak m=e^{-f}d\mu_{g|_\Sigma}\) and set
\[
\mathcal S_{N-1}^{\Sigma}=\mathcal S_{N-1}(g|_\Sigma,f|_\Sigma).
\]

\begin{lemma}\label{lem:stability}
On \(\Sigma_{\mathrm{reg}}\),
\begin{equation}\label{eq:weightedmean}
H-\nu f=\Phi
\end{equation}
and, for every \(v\in C_c^\infty(\Sigma_{\mathrm{reg}})\),
\begin{equation}\label{eq:stability}
\int_\Sigma\left(2|\nabla v|^2+\mathcal S_{N-1}^{\Sigma}v^2\right) d\mathfrak m\ge\int_\Sigma\left(\mathcal S_N(g,f)+\frac{N}{N-1}\Phi^2+2\nu\Phi\right)v^2 d\mathfrak m.
\end{equation}
\end{lemma}

\begin{proof}
The first and second variations of \eqref{eq:functional} give \eqref{eq:weightedmean} and
\begin{equation}\label{eq:secondvariation}
0\le\int_\Sigma\left[|\nabla v|^2-\left(\operatorname{Ric}_g(\nu,\nu)+\nabla^2f(\nu,\nu)+|A|^2+\nu\Phi\right)v^2\right] d\mathfrak m.
\end{equation}
The Gauss equation and \eqref{eq:weightedmean} give
\begin{equation}\label{eq:gauss}
\mathcal S_{N-1}^{\Sigma} =\mathcal S_N(g,f)-2\left(\operatorname{Ric}_g(\nu,\nu)+\nabla^2f(\nu,\nu)\right)+\Phi^2+\frac{(\nu f)^2}{N-m}-|A|^2.
\end{equation}
By \(|A|^2\ge H^2/(m-1)\) and \eqref{eq:weightedmean},
\begin{equation}\label{eq:trace}
|A|^2+\frac{(\nu f)^2}{N-m}-\frac{\Phi^2}{N-1} \ge\frac{N-1}{(m-1)(N-m)}\left(\nu f+\frac{N-m}{N-1}\Phi\right)^2 \ge0.
\end{equation}
Equations~\eqref{eq:secondvariation}--\eqref{eq:trace} give \eqref{eq:stability}.
\end{proof}

We obtain a pointwise lower bound for \(\mathcal S_N\) by changing the weight.

\begin{lemma}\label{lem:weightchange}
Let \((Y^{m-1},h,e^{-f}d\mu_h)\) be a smooth weighted Riemannian manifold, \(V\subset Y\) open, and \(K\Subset V\). Suppose \(w\in C^\infty(V)\) satisfies
\begin{equation}\label{eq:spectral}
\int_V\left(2|\nabla v|^2+(\mathcal S_{N-1}(h,f)-w)v^2\right) d\mathfrak m\ge0
\end{equation}
for every \(v\in C_c^\infty(V)\). There is \(\widehat f\in C^\infty(Y)\), equal to \(f\) outside \(V\), such that
\begin{equation}\label{eq:weightinequality}
\mathcal S_N(h,\widehat f)\ge w\quad\text{on }K.
\end{equation}
\end{lemma}

\begin{proof}
Choose an open set \(U\) with smooth boundary and \(K\Subset U\Subset V\). Set \(L=-2\Delta_f+\mathcal S_{N-1}(h,f)-w\). Let \(\psi>0\) be a first Dirichlet eigenfunction of \(L\) on each component of \(U\). Then \eqref{eq:spectral} gives \(L\psi\ge0\). Choose \(\widehat f\in C^\infty(Y)\) equal to \(f-\log\psi\) on a neighborhood of \(K\) and to \(f\) outside \(U\). On \(K\),
\begin{equation}\label{eq:weightsquare}
\mathcal S_N(h,\widehat f) =w+\frac{L\psi}{\psi} +\frac{N-m}{N-m+1}\left|d\log\psi+\frac{df}{N-m}\right|_h^2 \ge w.\qedhere
\end{equation}
\end{proof}

\section{Hardy inequalities near the singular set}\label{sec:hardy}

Let \(\Omega\) be a local minimizer of \eqref{eq:functional} in \(M^m\), with \(m\ge8\), and set \(\Sigma=\operatorname{spt}|D\chi_\Omega|\) and \(Z=\Sigma_{\mathrm{sing}}\). Then \(\dim_H Z\le m-8\).

We use \(B_r(x)\), \(\operatorname{dist}\) and \(\nabla_0\) for Euclidean balls, distance and connection.

The proof of the following lemma uses the stratification argument; see \cite[Sections~2--4]{WhiteStratification}.
\begin{lemma}\label{lem:flatpieces}
For every open set \(O\subset M\) and \(\eta>0\), \(Z\cap O\) is a countable union of compact sets \(F\). Each \(F\) lies in a coordinate chart with \(|g_{ij}-\delta_{ij}|<\eta\), and there is \(r_F>0\) such that
\begin{equation}\label{eq:flatness}
F\cap B_r(x)\subset\{z:\operatorname{dist}(z,L_{x,r})<\eta r\}, \qquad x\in F,\quad 0<r<r_F,
\end{equation}
where \(L_{x,r}\) is an affine \((m-8)\)-plane through \(x\).
\end{lemma}

\begin{proof}
It suffices to cover \(Z\cap Q\), where \(Q\Subset U\Subset O\) is compact and \(U\) is a coordinate chart with \(|g_{ij}-\delta_{ij}|<\eta\). Let \(\Theta(x)\) denote the density of \(\Sigma\) at \(x\). For \(x\in Z\cap Q\), there are \(d_x,r_x>0\), with \(B_{r_x}(x)\subset U\), such that, for every \(0<r<r_x\), some \((m-8)\)-plane \(L\) through \(x\) satisfies
\[
Z\cap B_r(x)\cap\{\Theta\ge\Theta(x)-d_x\} \subset\{z:\operatorname{dist}(z,L)<\eta r/4\}.
\]
If no such \(d_x,r_x\) exist, there are radii \(r_i\downarrow0\) for which the inclusion fails with \(d_x=1/i\). After passing to a subsequence, \(r_i^{-1}(\Sigma-x)\) converges locally as varifolds to a tangent cone \(T\) in the metric \(g(x)\). The monotonicity formula and its equality case give \(\Theta_T(z)\le\Theta_T(0)\) and
\[
\{z:\Theta_T(z)=\Theta_T(0)\} =\{z:T+z=T\}\subset\operatorname{sing}T.
\]
Since \(T\) is a cone, \(\{z:T+z=T\}\) is a linear subspace and has dimension at most \(m-8\).

Choose an \((m-8)\)-plane \(L\) containing \(\{z:T+z=T\}\) and \(y_i\in Z\cap B_{r_i}(x)\) with
\[
\operatorname{dist}\left((y_i-x)/r_i,L\right)\ge\eta/4,\qquad \Theta(y_i)\ge\Theta(x)-1/i.
\]
A limit point \(z\) of \((y_i-x)/r_i\) lies outside \(L\). Upper semicontinuity of density gives
\[
\Theta_T(z)\ge\Theta(x)=\Theta_T(0).
\]
Hence \(T+z=T\), contradicting \(z\notin L\).

For \(p\in\mathbb N\) and \(k\in\mathbb Z\), set
\[
A=\{x\in Z\cap Q:r_x,d_x\ge1/p,\ p\Theta(x)\in[k,k+1)\}, \qquad F=\overline A.
\]
Fix \(x\in F\) and \(0<r<1/(2p)\). Choose \(x_0\in A\) with \(|x-x_0|<\frac r2\min\{1,\eta\}\). Since \(2r<r_{x_0}\) and \(\Theta(y)>\Theta(x_0)-d_{x_0}\) for \(y\in A\), there is an \((m-8)\)-plane \(L\ni x_0\) with
\[
F\cap B_r(x)\subset\overline{A\cap B_{2r}(x_0)} \subset\{z:\operatorname{dist}(z,L)\le\eta r/2\}.
\]
Thus \eqref{eq:flatness} holds with \(L_{x,r}=L+x-x_0\) and \(r_F=1/(2p)\). The sets \(\overline A\), for all \(p,k\), cover \(Z\cap Q\).
\end{proof}

Fix a compact set \(F\subset Z\) as in Lemma~\ref{lem:flatpieces} and put \(\delta=\operatorname{dist}_g(\cdot,F)\). Write \(\|\cdot\|_2=\|\cdot\|_{L^2(\Sigma,d\mathfrak m)}\). In the next two lemmas, we prove Hardy inequalities with weight \(\delta^{-2}\). See \cite{Carron,Lehrback} for related estimates.

For \(w\in C^2(M)\), \eqref{eq:weightedmean} gives
\begin{equation}\label{eq:restrictedlaplacian}
\Delta_f^\Sigma w =\operatorname{tr}_{T\Sigma}^g\nabla_g^2w -\Phi\nu w-\langle df,dw\rangle_g.
\end{equation}
For \(0<U\in C^2(\Sigma_{\mathrm{reg}})\) and \(v\in C_c^\infty(\Sigma_{\mathrm{reg}})\),
\begin{equation}\label{eq:hardyidentity}
\int-\frac{\Delta_f^\Sigma U}{U}v^2 d\mathfrak m =\|\nabla v\|_2^2-\int U^2|\nabla(v/U)|^2 d\mathfrak m \le\|\nabla v\|_2^2.
\end{equation}

\begin{lemma}\label{lem:roughhardy}
There are \(\eta_0,c_0>0\), depending only on \(m\), such that for \(\eta\le\eta_0\),
\begin{equation}\label{eq:roughhardy}
c_0\int_\Sigma\delta^{-2}v^2 d\mathfrak m \le \|\nabla v\|_2^2+C_F\|v\|_2^2
\end{equation}
for every \(v\in C_c^\infty(\Sigma_{\mathrm{reg}})\), where \(C_F>0\) may depend on \(F\).
\end{lemma}

\begin{proof}
Fix \(\tau\in(0,1)\) small in terms of \(m\) and take \(\eta_0\le\tau/4\). By \eqref{eq:flatness}, for \(x\in F\) and \(0<r<r_F\), the set \(F\cap B_r(x)\) can be covered by at most \(C_m\tau^{-(m-8)}\le\tau^{-(m-7)}\) balls of radius \(\tau r\) centered in \(F\). After decreasing \(r_F\), iteration gives, for \(\delta(x)<r_F\) and \(0<t\le r<r_F\), a cover of \(F\cap B_{3r}(x)\) by at most \(C_m(r/t)^{m-7}\) balls of radius \(t\). Hence
\[
|B_r(x)\cap\{\delta<t\}| \le C_mr^{m-7}t^7,
\]
and
\begin{equation}\label{eq:distanceintegral}
\int_{B_r(x)}\delta(y)^{-6} dy \le C_m r^{m-6}.
\end{equation}

Set \(K(z)=|z|^{4-m}\). Choose a smooth function \(0\le\zeta\le1\), compactly supported in the chart and equal to one on \(\{\delta<2r_F\}\), and set
\[
U(x)=1+\zeta(x)\int\zeta(y)\delta(y)^{-6}K(x-y) dy.
\]
Set \(U=1\) outside the chart. Then \(U\in C^2(M\setminus F)\). For \(0<\delta(x)<r_F\) and \(j=0,1\), the inequality \(\delta(y)\ge\delta(x)/2\) on \(B_{\delta(x)/4}(x)\) gives
\[
\int_{B_{\delta(x)/4}(x)}\delta(y)^{-6}|\nabla_0^jK(x-y)|dy\le C_m\delta(x)^{-2-j}.
\]
For \(\delta(x)/4\le R<r_F/2\), \eqref{eq:distanceintegral} gives
\[
\int_{\{R<|x-y|<2R\}}\delta(y)^{-6}|\nabla_0^jK(x-y)|dy\le C_mR^{4-m-j}R^{m-6}=C_mR^{-2-j}.
\]
Thus, after decreasing \(r_F\),
\[
|\nabla_0^jU(x)|\le C_m\delta(x)^{-2-j},\qquad j=0,1.
\]
For every unit vector \(\xi\),
\[
-\operatorname{tr}_{\xi^\perp}\nabla_0^2K(z) =(m-4)|z|^{2-m} \left(1+(m-2)\frac{\langle z,\xi\rangle^2}{|z|^2}\right)\ge(m-4)|z|^{2-m}.
\]
Since \(|\nabla_0^2K(z)|\le C_m|z|^{2-m}\), choosing \(\eta_0\) small gives
\[
-\operatorname{tr}_{T_x\Sigma}^g\nabla_0^2U(x) \ge c_m\delta(x)^{-6}\int_{B_{\delta(x)/4}(x)}|x-y|^{2-m} dy \ge c_m\delta(x)^{-4}.
\]
By \eqref{eq:restrictedlaplacian},
\[
-\Delta_f^\Sigma U \ge c_m\delta^{-4}-C_F\delta^{-3} \quad \text{for }0<\delta<r_F,\quad -\frac{\Delta_f^\Sigma U}{U} \ge c_0\delta^{-2}-C_F \quad \text{on }\Sigma_{\mathrm{reg}}.
\]
Equation~\eqref{eq:hardyidentity} gives \eqref{eq:roughhardy}.
\end{proof}

\begin{lemma}\label{lem:annularhardy}
Fix \(A_0>1\) and \(0<h<(7-2)^2/4\). If \(\eta\) is sufficiently small in terms of \(m,A_0,h\), then for all small \(\varepsilon>0\),
\begin{equation}\label{eq:annularhardy}
h\int_{\{\varepsilon/2<\delta<A_0\varepsilon\}} \delta^{-2}v^2 d\mathfrak m \le \|\nabla v\|_2^2+C_F\|v\|_2^2
\end{equation}
for every \(v\in C_c^\infty(\Sigma_{\mathrm{reg}})\), where \(C_F>0\) is independent of \(\varepsilon\).
\end{lemma}

\begin{proof}
Fix \(t\in(0,\eta_0)\) and choose \(\Lambda\) with \(\log(t\Lambda/A_0)\ge1/t\). Assume \(\eta\le t/\Lambda\). Put \(r=\Lambda\varepsilon<r_F/4\).

For \(x_0\in F\), let \(L=L_{x_0,4r}\) in \eqref{eq:flatness} and set
\[
z^2=\operatorname{dist}(\cdot,L)^2+t^2\varepsilon^2, \qquad b=\frac{7-2}{2}.
\]
Since \(\dim L^\perp=8\), every hyperplane \(P\) satisfies
\[
\operatorname{tr}_P\nabla_0^2(z^2) =2\operatorname{tr}_P\operatorname{proj}_{L^\perp}\ge14.
\]
On \(B_{3r}(x_0)\), \eqref{eq:restrictedlaplacian} gives
\[
\Delta_f^\Sigma(z^2)\ge14-C_mt-C_Fr,\qquad |\nabla_\Sigma(z^2)|^2\le(4+C_mt)z^2.
\]
For \(C_Fr\le t\),
\[
\begin{aligned}
-z^{b+2}\Delta_f^\Sigma(z^{-b}) &=\frac b2\Delta_f^\Sigma(z^2) -\frac{b(b+2)}4\frac{|\nabla_\Sigma(z^2)|^2}{z^2}\\
&\ge b(7-2-b)-C_mt =\frac{(7-2)^2}{4}-C_mt>0.
\end{aligned}
\]

For \(x\in B_{3r}(x_0)\) with \(\varepsilon/2<\delta(x)<A_0\varepsilon\), choose \(y\in F\) with \(\operatorname{dist}_g(x,y)=\delta(x)\). Since \(y\in B_{4r}(x_0)\), \eqref{eq:flatness} gives
\[
z(x)\le|x-y|+4\eta r+t\varepsilon \le(1+C_mt)\delta(x).
\]

Choose a Lipschitz function \(0\le\chi=\chi(\delta)\le1\) with \(\chi=1\) for \(\delta\le A_0\varepsilon\), \(\chi=0\) for \(\delta\ge tr\), and \(|d\chi|\le t/\delta\). Choose finitely many points \(x_a\in F\) and smooth functions \(\varphi_a\) supported in \(B_{3r}(x_a)\) with \(\sum_a\varphi_a^2=1\) on \(\{\delta\le tr\}\) and \(\sum_a|d\varphi_a|^2\le C_mr^{-2}\). Using \(|d\chi|\le t/\delta\), \(r^{-2}\le t^2\delta^{-2}\) on \(\{0<\delta<tr\}\), and \eqref{eq:roughhardy},
\begin{align}
\sum_a\|\nabla(\chi\varphi_av)\|_2^2 &=\|\nabla(\chi v)\|_2^2+\int\chi^2\sum_a|\nabla\varphi_a|^2v^2 d\mathfrak m\notag\\
&\le(1+t)\|\nabla v\|_2^2+C_mt\int_{\{\delta<tr\}}\delta^{-2}v^2 d\mathfrak m\notag\\
&\le(1+C_mt)\|\nabla v\|_2^2+C_F\|v\|_2^2.\label{eq:localize}
\end{align}
Apply \eqref{eq:hardyidentity} to \(\chi\varphi_av\) with \(U=z^{-b}\) and \(x_0=x_a\). Summing over \(a\) and using \eqref{eq:localize} gives
\[
\left(\frac{25}{4}-C_mt\right) \int_{\{\varepsilon/2<\delta<A_0\varepsilon\}}\delta^{-2}v^2 d\mathfrak m \le(1+C_mt)\|\nabla v\|_2^2+C_F\|v\|_2^2.
\]
Choosing \(t\) so that \(25/4-C_mt>h(1+C_mt)\) proves \eqref{eq:annularhardy}.
\end{proof}

\begin{proposition}\label{prop:blowup-functions}
Assume
\begin{equation}\label{eq:N-range}
\frac{N-1}{N}<\frac{(7-2)^2}{4}.
\end{equation}
If \(\eta\) is sufficiently small in terms of \(m,N\), then, for some \(0<\rho<1\) and all small \(\varepsilon>0\), there are \(u_\varepsilon:M\to[-\infty,\infty]\) and \(0\le \sigma_\varepsilon\in C_c^\infty(\{0<\delta<\varepsilon\})\) satisfying
\begin{gather}
\begin{gathered}
F\subset\operatorname{int}\{u_\varepsilon=+\infty\},\qquad u_\varepsilon=-\infty\ \text{on }\{\delta\ge\varepsilon\},\\
2|du_\varepsilon|\le \frac{N}{N-1}u_\varepsilon^2+\sigma_\varepsilon \quad\text{a.e. where \(u_\varepsilon\) is finite},
\end{gathered}\label{eq:blowup-conditions}\\
\int \sigma_\varepsilon v^2 d\mathfrak m \le\rho \|\nabla v\|_2^2+C_F\|v\|_2^2\label{eq:blowup-error}
\end{gather}
for every \(v\in C_c^\infty(\Sigma_{\mathrm{reg}})\), where \(C_F>0\) is independent of \(\varepsilon\).
\end{proposition}

\begin{proof}
By \eqref{eq:N-range}, choose \(\beta,h,\rho>0\) with
\[
\frac{N-1}{N}(1+\beta^2)<\rho h<h<\frac{(7-2)^2}{4}.
\]
Put \(r=(\varepsilon/2)e^{-2\pi/\beta}\). The Riccati equation
\[
-2q'(t)=\frac{N}{N-1}q(t)^2+\frac{(N-1)(1+\beta^2)}{Nt^2}
\]
has the solution
\[
q(t)=\frac{N-1}{Nt}\left(1+\beta\cot\left(\frac\beta2\log(t/r)\right)\right),\qquad r<t<\varepsilon/2.
\]
See also \cite[Lemma~3.4]{Chen}. Set
\[
u_\varepsilon=\begin{cases}
+\infty,&\delta\le r,\\
q(\delta),&r<\delta<\varepsilon/2,\\
-\infty,&\delta\ge\varepsilon/2.
\end{cases}
\]
Since \(q'<0\) and \(|d\delta|\le1\),
\[
2|du_\varepsilon|\le\frac{N}{N-1}u_\varepsilon^2+\frac{(N-1)(1+\beta^2)}{N\delta^2}\quad\text{a.e. on }\{r<\delta<\varepsilon/2\}.
\]
Choose \(0\le \sigma_\varepsilon\in C_c^\infty(\{r/2<\delta<\varepsilon\})\) with
\[
\delta^2\sigma_\varepsilon\le\rho h\quad\text{on }M,\qquad \delta^2\sigma_\varepsilon\ge\frac{N-1}{N}(1+\beta^2) \quad\text{on }\{r\le\delta\le\varepsilon/2\}.
\]
Applying Lemma~\ref{lem:annularhardy} with \(r\) in place of \(\varepsilon\) and \(A_0=\varepsilon/r=2e^{2\pi/\beta}\) gives \eqref{eq:blowup-error}.
\end{proof}

\section{Dimension descent}

\begin{definition}\label{def:datum}
Let \(m\ge1\), \(N<0\) or \(m<N\le\infty\), and \(s>0\). An \(m\)-dimensional \emph{descent datum} for \(N,s\) consists of a compact smooth oriented weighted \(m\)-manifold \((X_m,g,f_m)\) and functions \(\Phi_1,\ldots,\Phi_m\) defined on a neighborhood of \(X_m\). The functions \(\Phi_j\) are continuous with values in \([-\infty,\infty]\) and locally Lipschitz wherever finite. There exist smooth positive functions \(\sigma_j\) on a neighborhood of \(X_m\) and relatively open sets \(U_j^\pm\subset X_m\), \(1\le j\le m\), such that
\begin{equation}\label{eq:Phi-conditions}
2|d\Phi_j|\le \frac{N}{N-1}\Phi_j^2+\sigma_j \quad\text{a.e. on \(X_m\cap\{|\Phi_j|<\infty\}\)},\qquad \overline{U_j^\pm}\subset\operatorname{int}\{\Phi_j=\pm\infty\},
\end{equation}
and, with \(K_m=X_m\setminus\bigcup_{j=1}^m(U_j^+\cup U_j^-)\),
\begin{equation}\label{eq:datum-conditions}
K_m\Subset\operatorname{int} X_m,\qquad \mathcal S_N(g,f_m)\ge s+\sum_{j=1}^m \sigma_j\quad\text{on }K_m,\qquad \deg(\Phi_1,\ldots,\Phi_m)\ne0.
\end{equation}
\end{definition}

Fix such a datum, choose \(\sigma_j,U_j^\pm\) as above, and let \(0<\varepsilon<s\). Extend \(g,f_m\) smoothly to a closed \(m\)-manifold \(M\) containing \(X_m\), and write \(d\mathfrak m=e^{-f_m}d\mu_g\).
\begin{lemma}\label{lem:minimizer}
There are \(\Phi\in C^\infty(M)\) and a minimizer \(\Omega\) of \(\mathcal J(\cdot;M)\) such that
\begin{gather}
\mathcal S_N(g,f_m)+\frac{N}{N-1}\Phi^2-2|d\Phi|>s+\sum_{j<m}\sigma_j-\frac\varepsilon2\quad\text{on }K_m,\label{eq:minimizer-curvature}\\
\overline{U_m^+}\subset\operatorname{int}\Omega,\qquad \overline{U_m^-}\subset\operatorname{int}(M\setminus\Omega).\label{eq:minimizer-inclusions}
\end{gather}
\end{lemma}

\begin{proof}
Choose open sets \(D^\pm\subset M\) with smooth boundary and disjoint closures such that
\[
\overline{U_m^\pm}\subset D^\pm,\qquad \overline{D^\pm}\cap X_m\subset\operatorname{int}\{\Phi_m=\pm\infty\}.
\]
For a vector field \(\xi\), set \(\operatorname{div}_{f_m}\xi=\operatorname{div}_g\xi-df_m(\xi)\). Let \(\xi^\pm\) be smooth vector fields with \(|\xi^\pm|\le1\), equal to the outward unit normal on \(\partial D^\pm\), and fix \(L>1+\max_\pm\|\operatorname{div}_{f_m}\xi^\pm\|_\infty\).

By \eqref{eq:Phi-conditions}, smooth the truncation \(\max\{-L,\min\{\Phi_m,L\}\}\) to \(\Phi\in C^\infty(M)\) with
\begin{equation}\label{eq:smooth-Phi}
\Phi=\pm L\text{ on }D^\pm,\qquad 2|d\Phi|<\frac{N}{N-1}\Phi^2+\sigma_m+\varepsilon/2\text{ on }K_m.
\end{equation}
Equation~\eqref{eq:datum-conditions} then gives \eqref{eq:minimizer-curvature}.

Let \(\Omega\) minimize \(\mathcal J(\cdot;M)\) and write \(P_{f_m}(E)=\int_{\partial^*E}d\mathfrak m\). Then
\[
P_{f_m}(\Omega\cup D^+)-P_{f_m}(\Omega)\le P_{f_m}(D^+)-P_{f_m}(\Omega\cap D^+)\le\int_{D^+\setminus\Omega}\operatorname{div}_{f_m}\xi^+ d\mathfrak m.
\]
Hence
\[
0\le\mathcal J(\Omega\cup D^+;M)-\mathcal J(\Omega;M)\le\int_{D^+\setminus\Omega}\left(\operatorname{div}_{f_m}\xi^+-L\right)d\mathfrak m\le-\mathfrak m(D^+\setminus\Omega).
\]
The same estimate for \(M\setminus\Omega\) and \(-\Phi\) gives \(\mathfrak m(D^-\cap\Omega)=0\), proving \eqref{eq:minimizer-inclusions}.
\end{proof}

Assume \(m\ge2\), and assume \eqref{eq:N-range} if \(m\ge8\). Fix \(\Phi,\Omega\) as in Lemma~\ref{lem:minimizer} and set \(\Sigma=\operatorname{spt}|D\chi_\Omega|\), \(Z=\Sigma_{\mathrm{sing}}\). Choose an open set \(O\) with \(K_m\subset O\Subset\operatorname{int} X_m\) on which \eqref{eq:minimizer-curvature} holds. Put \(\mathcal W=W_0^{1,2}(\Sigma_{\mathrm{reg}}\cap O,d\mathfrak m)\). For \(v\in C_c^\infty(\Sigma_{\mathrm{reg}}\cap O)\), set
\[
\mathcal Q(v)=\int_{\Sigma\cap O} \left[2|\nabla v|^2+ \left(\mathcal S_{N-1}^{\Sigma}-s+\varepsilon-\sum_{j<m}\sigma_j\right)v^2\right] d\mathfrak m.
\]
By \eqref{eq:secondvariation} and \eqref{eq:gauss}, respectively,
\[
\int_\Sigma|A|^2v^2d\mathfrak m\le\|\nabla v\|_2^2+C\|v\|_2^2,\qquad \mathcal S_{N-1}^{\Sigma}+|A|^2\in L^\infty(\Sigma).
\]
Hence \(\mathcal Q\) extends continuously to \(\mathcal W\). Together with \eqref{eq:stability} and \eqref{eq:minimizer-curvature}, these estimates give, for some \(C_0\ge0\),
\begin{equation}\label{eq:formbounds}
\mathcal Q(v)\ge\frac\varepsilon2\|v\|_2^2,\qquad \|\nabla v\|_2^2-C_0\|v\|_2^2 \le\mathcal Q(v) \le2\|\nabla v\|_2^2+C_0\|v\|_2^2.
\end{equation}

\begin{lemma}\label{lem:error-absorption}
Assume \(m\ge8\). Let \(\sigma_\tau\) and \(0<\rho<1\) be given by Proposition~\ref{prop:blowup-functions} for a fixed \(F\subset Z\). Suppose \(V\in L^\infty(\Sigma,d\mathfrak m)\) and \(a<1\) satisfy
\[
\int Vv^2d\mathfrak m\le a\mathcal Q(v)\qquad\text{for every }v\in\mathcal W.
\]
For every \(\max\{a,\rho\}<b<1\), all sufficiently small \(\tau>0\), and every \(v\in\mathcal W\),
\begin{equation}\label{eq:error-absorption}
\int(V+\sigma_\tau)v^2d\mathfrak m\le b\mathcal Q(v).
\end{equation}
\end{lemma}

\begin{proof}
By \eqref{eq:formbounds}, \(\mathcal Q^{1/2}\) is a norm equivalent to the \(W^{1,2}\)-norm on \(\mathcal W\). The mean curvature of \(\Sigma\) is bounded, so the Michael--Simon Sobolev inequality~\cite{MichaelSimon} gives
\[
\|v\|_{L^{2(m-1)/(m-3)}(d\mathfrak m)}^2\le C\left(\|\nabla v\|_2^2+\|v\|_2^2\right)\le C\mathcal Q(v).
\]
For \(E_r=\Sigma\cap\{\operatorname{dist}_g(\cdot,Z)<r\}\), H\"older's inequality gives
\[
\sup_{\mathcal Q(v)\le1}\int_{E_r}v^2d\mathfrak m\le C\mathfrak m(E_r)^{2/(m-1)}\longrightarrow0\quad\text{as }r\downarrow0,
\]
since \(\mathfrak m(Z)=0\). Together with local Rellich compactness for zero extensions to \(\Sigma_{\mathrm{reg}}\), this proves that \(\mathcal W\) embeds compactly in \(L^2(d\mathfrak m)\).

For fixed \(v\in\mathcal W\) and \(w\in C_c^\infty(\Sigma_{\mathrm{reg}}\cap O)\), \(\sigma_\tau w=0\) for small \(\tau\), so \eqref{eq:blowup-error} gives
\[
\limsup_{\tau\downarrow0}\int\sigma_\tau v^2d\mathfrak m\le\rho\|\nabla(v-w)\|_2^2+C_F\|v-w\|_2^2.
\]
Taking \(w\to v\) in \(\mathcal W\) shows that \(\int\sigma_\tau v^2d\mathfrak m\to0\).

If \eqref{eq:error-absorption} fails for arbitrarily small \(\tau\), choose \(\tau_i\downarrow0\) and \(v_i\in\mathcal W\) with \(\mathcal Q(v_i)=1\) and \(\int(V+\sigma_{\tau_i})v_i^2d\mathfrak m>b\). After passage to a subsequence, \(v_i\rightharpoonup v\) in \(\mathcal W\) and \(v_i\to v\) in \(L^2(d\mathfrak m)\). Weak convergence gives
\[
\mathcal Q(v_i-v)=1-\mathcal Q(v)+o(1).
\]
By \eqref{eq:formbounds} and \eqref{eq:blowup-error},
\[
\begin{aligned}
\int\sigma_{\tau_i}(v_i-v)^2d\mathfrak m&\le\rho\mathcal Q(v_i-v)+(C_F+\rho C_0)\|v_i-v\|_2^2\\
&=\rho\left(1-\mathcal Q(v)\right)+o(1).
\end{aligned}
\]
Thus,
\[
\left|\int\sigma_{\tau_i}\left(v_i^2-(v_i-v)^2\right)d\mathfrak m\right|\le C\left(\int\sigma_{\tau_i}v^2d\mathfrak m\right)^{1/2}+\int\sigma_{\tau_i}v^2d\mathfrak m\longrightarrow0.
\]
Since \(0\le\mathcal Q(v)\le1\),
\begin{align*}
b&\le\limsup_i\int(V+\sigma_{\tau_i})v_i^2d\mathfrak m\\
&\le a\mathcal Q(v)+\rho\left(1-\mathcal Q(v)\right)\le\max\{a,\rho\}<b.\qedhere
\end{align*}
\end{proof}

\begin{proposition}\label{prop:Phi-modification}
There are functions \(\widehat\Phi_j,\widehat\sigma_j\) and open sets \(\widehat U_j^\pm\supset U_j^\pm\), \(j<m\), satisfying \eqref{eq:Phi-conditions}, with \(0\le \widehat\sigma_j-\sigma_j\in C_c^\infty(O)\) and
\begin{equation}\label{eq:singular-cover}
Z\cap X_m\subset\bigcup_{j<m}(\widehat U_j^+\cup \widehat U_j^-).
\end{equation}
For every \(0\ne v\in\mathcal W\),
\begin{equation}\label{eq:totalerror}
\sum_{j<m}\int_\Sigma(\widehat\sigma_j-\sigma_j)v^2 d\mathfrak m <\mathcal Q(v).
\end{equation}
\end{proposition}

\begin{proof}
By \eqref{eq:minimizer-inclusions},
\[
K:=Z\cap X_m\setminus\bigcup_{j<m}(U_j^+\cup U_j^-)=Z\cap K_m\subset O.
\]
We may assume \(K\ne\varnothing\), so \(m\ge8\).

Since \(\dim_H K<m-1\) and each \(\Phi_j\) is locally Lipschitz where finite, choose \(t\in\mathbb R^{m-1}\setminus(\Phi_j)_{j<m}(K)\). Set \(O_j^\pm=\{\pm(\Phi_j-t_j)>0\}\). Then
\[
K\subset\bigcup_{j<m}(O_j^+\cup O_j^-).
\]
Apply Lemma~\ref{lem:flatpieces} in \(O\cap O_j^\pm\), with \(\eta\) small in terms of \(m,N\), to cover \(K\) by compact sets \(F_a\).

Set \(b_0=0\). For each \(a\), apply Proposition~\ref{prop:blowup-functions} to \(F_a\) and choose \(b_a\) with \(\max\{b_{a-1},\rho\}<b_a<1\). For sufficiently small \(\tau\), Lemma~\ref{lem:error-absorption} with \(V=\sum_{k<a}\widetilde \sigma_k\) gives \(u_a=u_\tau\) and \(\widetilde \sigma_a=\sigma_\tau\in C_c^\infty(O)\) with \(\{u_a>-\infty\}\Subset O\cap O_j^\pm\) for some \(O_j^\pm\) and
\[
\sum_{k\le a}\int_\Sigma\widetilde \sigma_k v^2 d\mathfrak m\le b_a\mathcal Q(v)
\]
for every \(v\in\mathcal W\).

By compactness of \(K\), choose \(R\) such that
\[
K\subset\bigcup_{a\le R}\operatorname{int}\{u_a=+\infty\}.
\]
Partition \(\{1,\ldots,R\}=\bigsqcup_{j<m}(I_j^+\sqcup I_j^-)\) so that, for \(a\in I_j^\pm\),
\[
\{u_a>-\infty\}\subset O_j^\pm.
\]
Define
\begin{equation}\label{eq:modification}
\widehat\Phi_j=\min\left\{ \max\{\Phi_j,u_a:a\in I_j^+\},-u_a:a\in I_j^-\right\},\qquad \widehat\sigma_j=\sigma_j+\sum_{a\in I_j^+\cup I_j^-}\widetilde \sigma_a.
\end{equation}
At almost every point where \(\widehat\Phi_j\) is finite, its value and differential agree with those of \(\Phi_j\), some \(u_a\), or some \(-u_a\). By \eqref{eq:blowup-conditions}, the differential inequality in \eqref{eq:Phi-conditions} holds for \(\widehat\Phi_j\) and \(\widehat\sigma_j\). Since \(O_j^+\cap O_j^-=\varnothing\), choose open sets \(\widehat U_j^\pm\supset U_j^\pm\) with \(\overline{\widehat U_j^\pm}\subset\operatorname{int}\{\widehat\Phi_j=\pm\infty\}\) and satisfying \eqref{eq:singular-cover}.
\end{proof}

\begin{remark}\label{rem:sixdimensional}
The proof only uses \(\dim_H Z<m-1\) to choose \(t\). Using \(\dim_H Z\le m-8\), one can strengthen Proposition~\ref{prop:Phi-modification} with \(\widehat\Phi_j=\Phi_j\), \(\widehat\sigma_j=\sigma_j\) and \(\widehat U_j^\pm=U_j^\pm\) for \(j<m\) with \(j\le6\).
\end{remark}

\begin{proposition}\label{prop:descent}
There is a descent datum
\[
(X_{m-1},g|_{X_{m-1}},f_{m-1};\widehat\Phi_1,\ldots,\widehat\Phi_{m-1})
\]
for \(N,s-\varepsilon\), where \(X_{m-1}\subset X_m\) is a hypersurface with \(\partial X_{m-1}=X_{m-1}\cap\partial X_m\).
\end{proposition}

\begin{proof}
Choose \(\widehat\Phi_j,\widehat\sigma_j\) and \(\widehat U_j^\pm\) by Proposition~\ref{prop:Phi-modification}. Choose an open set \(V\subset\Sigma_{\mathrm{reg}}\) with
\[
K_{m-1}:=\Sigma\cap X_m\setminus \bigcup_{j<m}(\widehat U_j^+\cup \widehat U_j^-) \Subset V\Subset\Sigma_{\mathrm{reg}}\cap O.
\]
Let \(d\) be the signed distance to \(\Sigma\), positive in \(\Omega\). Choose \(h\in C^\infty(M,I)\) equal to \(d\) on a neighborhood of \(\overline V\) and to \(\pm1\) on neighborhoods of \(\overline{U_m^\pm}\), with
\[
hd>0\quad\text{on }X_m\setminus \left(\Sigma\cup\bigcup_{j<m}(\widehat U_j^+\cup \widehat U_j^-)\right).
\]
Choose \(h\) so that \(0\) is a regular value of both \(h|_{X_m}\) and \(h|_{\partial X_m}\). Set \(X_{m-1}=\{h=0\}\cap X_m\). Then
\[
X_{m-1}\setminus\bigcup_{j<m}(\widehat U_j^+\cup \widehat U_j^-) =K_{m-1}\Subset\operatorname{int} X_{m-1}.
\]
For every \(v\in C_c^\infty(V)\), \eqref{eq:totalerror} gives
\[
\int_V\left(2|\nabla v|^2+\left(\mathcal S_{N-1}^{\Sigma}-s+\varepsilon-\sum_{j<m}\widehat\sigma_j\right)v^2\right)d\mathfrak m\ge0.
\]
Since \(X_{m-1}\) agrees with \(\Sigma\) on \(V\), Lemma~\ref{lem:weightchange} with \(w=s-\varepsilon+\sum_{j<m}\widehat\sigma_j\) gives \(f_{m-1}\) with \(f_{m-1}-f_m|_{X_{m-1}}\in C_c^\infty(V)\) and
\begin{equation}\label{eq:nextcurvature}
\mathcal S_N(g|_{X_{m-1}},f_{m-1})\ge s-\varepsilon+\sum_{j<m}\widehat\sigma_j \quad\text{on }K_{m-1}.
\end{equation}

For \(x\in\operatorname{int} X_m\) with \(|\widehat\Phi_j(x)|<\infty\), consider a small geodesically convex ball about \(x\) on which \(\widehat\Phi_j\) is finite. By \eqref{eq:Phi-conditions}, the supremum of \(\frac{N}{N-1}\widehat\Phi_j^2+\widehat\sigma_j\) on this ball bounds twice the Lipschitz constant of \(\widehat\Phi_j\). Letting the radius tend to zero and using continuity gives
\[
2\limsup_{y\to x} \frac{|\widehat\Phi_j(y)-\widehat\Phi_j(x)|}{\operatorname{dist}_g(x,y)} \le \frac{N}{N-1}\widehat\Phi_j(x)^2+\widehat\sigma_j(x).
\]
Since \(\operatorname{dist}_g\le\operatorname{dist}_{g|_{X_{m-1}}}\), \(\widehat\Phi_j,\widehat\sigma_j\) satisfy the differential inequality in \eqref{eq:Phi-conditions} on \(X_{m-1}\).

The sets \(U_j^\pm\) cover \(\partial X_m\), with \(\Phi_j=\widehat\Phi_j=\pm\infty\) on \(U_j^\pm\) for \(j<m\) and \(\Phi_m=\pm\infty\), \(h=\pm1\) on \(U_m^\pm\). Orient \(X_{m-1}\) by \(TX_m=TX_{m-1}\oplus\mathbb R\nabla h\). By homotopy invariance and computation of the degree on \(X_{m-1}=h^{-1}(0)\cap X_m\),
\[
\deg\left((\Phi_j)_{j\le m}\right) =\deg\left((\widehat\Phi_j)_{j<m},h\right)=\deg\left((\widehat\Phi_j)_{j<m}|_{X_{m-1}}\right).\qedhere
\]
\end{proof}

\begin{proof}[Proof of Theorem~\ref{thm:cube}]
By \cite[Theorem~B]{BZ} with \(\lambda=1/N\), we may assume \(\frac{N-1}{N}<\frac{25(n-1)}{23n+26}<\frac{25}{4}\). Suppose \(\sigma>4\pi^2\frac{N-1}{N}\sum_i d_i^{-2}\). Choose \(0<\bar d_i<d_i\) so that
\[
\sigma_i=\frac{4\pi^2(N-1)}{N\bar d_i^2},\qquad \sum_i\sigma_i<\sigma.
\]
Put \(t_i=\operatorname{dist}_g(\cdot,B_i^-)-(d_i-\bar d_i)/2\) and
\[
\Phi_i=\begin{cases}
-\infty,&t_i\le0,\\[2pt]
\displaystyle-\frac{2\pi(N-1)}{N\bar d_i} \cot \left(\frac{\pi t_i}{\bar d_i}\right),&0<t_i<\bar d_i,\\[6pt]
+\infty,&t_i\ge\bar d_i.
\end{cases}
\]
Extend \(\Phi_i\) to a neighborhood of \(X\) with \(B_i^\pm\subset\operatorname{int}\{\Phi_i=\pm\infty\}\). The identity \(2 d\Phi_i/dt_i=\frac{N}{N-1}\Phi_i^2+\sigma_i\) and \(|dt_i|\le1\) give \eqref{eq:Phi-conditions}. Choose \(U_i^\pm\supset B_i^\pm\) with \(\overline{U_i^\pm}\subset\operatorname{int}\{\Phi_i=\pm\infty\}\). Since \(\Phi_i=\pm\infty\) and \(F_i=\pm1\) on \(B_i^\pm\), homotopy invariance gives
\[
\deg(\Phi_1,\ldots,\Phi_n)=\deg F\ne0.
\]

With \(X_n=X\), \(f_n=f\) and \(s=\sigma-\sum_i\sigma_i>0\), we obtain an \(n\)-dimensional descent datum for \(N,s\). Apply Proposition~\ref{prop:descent} \(n-1\) times with \(\varepsilon=s/n\) and fixed \(N\).

Let \(\Phi,\Omega\) be given by Lemma~\ref{lem:minimizer} on \(X_1\), with \(0<\varepsilon<s/n\). Since \(\deg\Phi_1\ne0\), some component of \(X_1\) has one endpoint in \(U_1^-\) and the other in \(U_1^+\). By \eqref{eq:minimizer-inclusions}, its interior contains a point \(x\in\partial^*\Omega\cap K_1\). At \(x\), the first and second variations give
\[
-\nu f_1=\Phi,\qquad \nabla^2f_1(\nu,\nu)+\nu\Phi\le0.
\]
Equation~\eqref{eq:minimizer-curvature} then gives
\[
\begin{aligned}
0&<\mathcal S_N(g,f_1)+\frac{N}{N-1}\Phi^2-2|d\Phi|\\
&\le\mathcal S_N(g,f_1)+\frac{N}{N-1}\Phi^2+2\nu\Phi=2\left(\nabla^2f_1(\nu,\nu)+\nu\Phi\right)\le0,
\end{aligned}
\]
a contradiction.

For sharpness, take \(0<a<\pi/2\) and \(L>0\), and set
\[
g=dt^2+(\cos t)^{4/N}|dx|^2,\qquad f=2\left(\frac nN-1\right)\log\cos t
\]
on \([-a,a]\times[-L,L]^{n-1}\). Then \(\mathcal S_N(g,f)=4(N-1)/N\) and \(d_1=2a\). For~\(i\ge2\),
\[
d_i\ge2L\min_{|t|\le a}(\cos t)^{2/N}.
\]
Letting \(L\to\infty\) and then \(a\uparrow\pi/2\) proves that the coefficient \((N-1)/N\) in \eqref{eq:cube} is optimal.
\end{proof}

\end{document}